\documentclass{article}

\usepackage{graphicx}          

\usepackage{color} 

 \newcommand\figref{Figure~\ref}

 \usepackage{comment,amsmath,gensymb,booktabs}
\usepackage{dsfont}\let\mathbb=\mathds

\usepackage{standalone}
\usepackage{amsthm}
\usepackage{mathtools,amsmath}
\usepackage{graphicx}
\usepackage{enumerate}
\usepackage{tikz}
\usepackage{color}
\usepackage{comment}
\usepackage{bibentry}

\usepackage{pgfplots}
\usetikzlibrary{pgfplots.groupplots}
\usetikzlibrary{fit}
\pgfplotsset{compat=1.9}
\usetikzlibrary{positioning}
\usepackage{graphicx}
\newcommand{\rev}{\textcolor{black}}
\newtheorem{defn}{Definition}
\newtheorem{lem}{Lemma}
\newtheorem{thm}{Theorem}
\newtheorem{rem}{Remark}
\newtheorem{ack}{Acknowledgements}

\begin{document}


\title{\rev{Closed-form H-infinity Optimal Control for a \\ Class of Infinite-Dimensional Systems}\\[15pt] \large Carolina Bergeling, Kirsten A. Morris$^*$ and Anders Rantzer \\[20pt] \small Department of Automatic Control, Lund University\\ Box 118, SE-221 00 Lund, Sweden. \\[5pt] 
$^*$Department of Applied Mathematics, University
of Waterloo\\ Ontario, Canada N2L 3G1. \\[10pt] 
Corresponding author: Carolina Bergeling\\ carolina.bergeling@control.lth.se, +46 46 222 87 87. } 
\date{}                                                
\maketitle

\section*{Abstract}                          
\rev{H-infinity optimal control and estimation are addressed for a class of systems governed by partial differential equations with bounded input and output operators.} Diffusion equations are an important example in this class. \rev{Explicit formulas for the} optimal state feedback controller as well as the optimal state estimator are given. Unlike traditional methods for H-infinity synthesis, no iteration is needed to obtain the optimal solution. Moreover, the optimal performance for both the state feedback and state estimation problems are  explicitly calculated. This is shown to be useful for problems of H-infinity optimal actuator and sensor location. Furthermore, the results can be used in testing and bench-marking of general purpose algorithms for H-infinity synthesis. \rev{The results also apply to finite-dimensional systems. }

\section{Introduction}

There are systems for which the variable of interest varies in both time and space. These systems are modelled by partial differential equations. Since the state evolves on an infinite-dimensional space, they are often referred to as infinite-dimensional systems. In controller design for infinite-dimensional systems, the synthesis problem is often approached by first approximating the partial differential equations by a system of ordinary differential equations. Then, the synthesis is performed based on this finite-dimensional approximation of the original system. Methods with this approach are called indirect, see for example \cite{Morris-DPSbook} for a review of indirect approaches to LQ and $H_{\infty}$ control.  On the contrary, direct methods perform synthesis directly on the infinite-dimensional system. Such methods are treated in \cite{curtain2012introduction},\cite{Ozbay-etal-book},\cite{Ke}.

The difficulty with indirect methods lies in ensuring that the designed controller works as predicted on the original infinite-dimensional system. The advantage of a direct approach is that since the controller is designed for the original model, actual performance is easier to determine.  However, the direct approach is in general difficult or impossible to use as it usually requires a closed-form description of the system, which is generally not  available. 

In this work, $H_{\infty}$ optimal control of infinite-dimensional systems is addressed. The aim of $H_{\infty}$ control is to stabilize a system as well as attenuate its response to worst-case disturbances. This is an alternative to for instance LQG, where the disturbances are assumed to be Gaussian white noise. For finite-dimensional systems, the $H_{\infty}$ optimal synthesis problem is generally treated by iteratively solving a series of algebraic Riccati equations (AREs), see e.g. \cite{Petersen1987} and \cite{doyle1989state}.  This approach can be computationally intensive and it is not uncommon for numerical issues to arise, especially  when attempting to compute a controller near optimal attenuation \cite{lanzon2008computing}. This is particularly true for large-scale systems such as high order approximations of infinite-dimensional systems \cite{kasinathan2014solution}. Furthermore, for infinite-dimensional systems, it needs to be verified that the controller synthesized for the finite-dimensional approximation is admissible for the original infinite-dimensional system, as previously explained. Sufficient conditions for this to be the case are described in \cite{ito1998approximation,Mo01}. 

This work presents a method for direct calculation of the optimal $H_\infty$ attenuation and the optimal control law for a class of infinite-dimensional systems. In fact, both the optimal attenuation and the optimal control law are stated explicitly on closed form. Approximation, if needed, is done on the explicitly stated control law. The method presented is applicable to a class of systems in which diffusion equations are an important example. 
The synthesis problems treated are state feedback and its dual, filtering. Moreover, the results are illustrated with several examples. Explicit expressions for the optimal state feedback and estimations laws are given for diffusion problems with constant as well as variable conductivity, in both one and two spatial dimensions. Similar results for a class of finite-dimensional systems were obtained in \cite{lidstrom2016optimal}.

It is rare to obtain the solution of an $H_{\infty}$ optimal control problem on closed form. This is true for both finite and infinite dimensional systems. In general, gamma-iteration, i.e. a bisection algorithm, is needed to find the optimal controller. The one-block problem, or one-block Nehari problem, is one exception where the optimal controller can be given in one shot, only dependent on the solutions of two Lyapunov equations \cite[Ch. 4]{ionescu2012robust}. The closed-form solution presented in this work is another example. However, the controller is given directly as well as is in no need of iteration to guarantee optimal performance. \rev{To the authors's knowledge, there are no existing results of this form besides that previously presented by the authors in \cite{lidstrom2016optimal} for finite-dimensional systems.} Moreover, the explicitly stated controller and estimator are shown to be very powerful tools. Both for their primary purpose of control and estimation, respectively, and as a means for benchmarking of general purpose algorithms for $H_{\infty}$ control, especially for systems of large order. Moreover, the results are suggested as tools for the problems of optimal actuator and sensor placement. 

The outline is as follows. Section~\ref{math} provides the mathematical framework. In Section~\ref{theorysection1} and \ref{theorysection2}, the closed-form $H_{\infty}$ optimal state feedback and estimation laws are stated. \rev{A preliminary version of the state feedback result was given in \cite{lidstrom2016h}. However, this work provides a completely new proof which illustrates the underlying properties of the class of systems that can be treated, namely that they are most vulnerable to constant disturbances. The state feedback and estimation results are illustrated with analytic examples in the sections thereafter.} Finally, Section~\ref{numericalsection} discusses how the results can be used for general purpose algorithms for $H_{\infty}$ control. 
 
\section{Notation and Definitions}
\label{math}

The set of real and complex numbers are denoted $\mathbb{R}$ and $\mathbb{C}$, respectively.  Furthermore, $\mathbb{R}_+$ is the set of nonnegative real numbers. The notation $\textrm{Re}(x)$, where $x \in \mathbb{C}$, denotes the real part of $x$. The scalar product is denoted $\langle \cdot , \cdot \rangle$ and a norm is denoted by $\|\cdot\|$, where a subscript is used when in need of clarification. For instance, $\|\cdot\|_2$ denotes the norm on $L_2[0,\infty)$.  The Laplace transform of a function $y$ is denoted $\hat{y}$. 

Only linear operators on separable Hilbert spaces will be considered. Given a linear operator $$T \, : D(T) \subset \mathcal{X} \rightarrow \mathcal{Y},$$
where $D(T)$ is the domain of $T$, its adjoint is denoted $T^*$ and its inverse, if it exists, is denoted $T^{-1}$. The set of bounded linear operators $T$ from $\mathcal{X}$ to $\mathcal{Y}$ is denoted $\mathcal{L}\,(\mathcal{X},\mathcal{Y})$, and ${\mathcal{L}\,(\mathcal{X}) \coloneqq \mathcal{L}\,(\mathcal{X},\mathcal{X})}$. The following definitions and results are helpful in the proceeding sections. For more details  consider \cite{lax}, \cite{luenberger1969} and \cite{curtain2012introduction}.

\begin{defn}[{\cite{curtain2012introduction}}]
A self-adjoint operator $A$ on the Hilbert space $\mathcal{Z}$ is nonnegative if ${\langle Az,z \rangle \geq 0} $ for all $ z \in D(A)$. The operator $A$ is positive if ${\langle Az,z \rangle > 0}$ for all nonzero $ z \in D(A)$. Furthermore, $A$ is {\em strictly positive (or coercive)} if there exists an $m > 0$ such that 
 \begin{equation*} \langle Az,z \rangle \geq m \|z\|^2 \text{ for all } z \in D(A). 
  \end{equation*}
The operator $A$ is (strictly) negative if $-A$ is (strictly) positive. 
\end{defn}

If $A: D(A) \subset \mathcal Z \to \mathcal Z$ is the infinitesimal generator of a strongly continuous semigroup $T(t)$ on $\mathcal{Z}$, then for all $z_0 \in D(A)$, the following differential equation on $\mathcal{Z}$
\begin{align} \label{cauchyproblem}
\frac{dz(t)}{dt} = Az(t), t\geq 0\quad z(0) = z_0,
\end{align}
has the unique solution $z(t) = T(t)z_0$, and this solution depends continuously on the initial condition. See \cite[p. 21, Th. 2.1.10]{curtain2012introduction}. Thus, the differential equation is well-posed. 

The following theorem is useful for determining well-posedness.
 \begin{lem}[{\cite[p. 33, Cor. 2.2.3]{curtain2012introduction}}] 
 \label{lem:1}\label{lem1}%
Sufficient conditions for a closed, densely defined operator $A$ on a Hilbert space to be the infinitesimal generator of a strongly continuous semigroup satisfying $\|T(t)\|\leq e^{\omega t}$ are 
\begin{align*}
\text{\normalfont Re}(\langle Az,z \rangle) &\leq \omega \|z\|^2 \quad \text{ for } z \in D(A), \\
\textrm{\normalfont Re}(\langle A^*z,z \rangle) &\leq \omega \|z\|^2 \quad \text{ for } z \in D(A^*).
\end{align*} 
\end{lem}
\begin{defn}[{\cite[p. 215, Def. 5.1.1]{curtain2012introduction}}]
A strongly continuous semigroup $T(t)$ on a Hilbert space $\mathcal{Z}$  is  \textit{exponentially stable} if there exist positive constants $M,\, \alpha$ such that ${\|T(t)\|\leq M e^{-\alpha t}}$ for $t \geq 0$. \label{def:exp}
\end{defn}
The following theorem will be useful for proving exponential stability. 

\begin{thm}[{\cite[p. 217, Th.~5.1.3]{curtain2012introduction}}] 
\label{lemLyap}%
Suppose that ${A: D(A)\subset \mathcal{Z} \rightarrow \mathcal{Z}}$ is the infinitesimal generator of a strongly continuous semigroup $T(t)$ on the Hilbert space $\mathcal{Z}$. Then $T(t)$ is exponentially stable if and only if there exists a positive operator $P\in \mathcal{L}(\mathcal Z)$ such that 
$$\langle Az,Pz \rangle +\langle Pz,Az \rangle \leq -\langle z,z\rangle$$ 
for all $z\in D(A)$.
\end{thm}

Consider the causal control system
\begin{align} \label{sys:AB}
\frac{dz(t)}{dt} = Az(t)+Bu(t), \quad z(0) = z_0, \quad t\geq 0,
\end{align}
where $A$ is  the infinitesimal generator of a strongly continuous semigroup $T(t)$ on $\mathcal{Z}$ while $B$ is a bounded operator from the input space $\mathcal{U}$ to $\mathcal{Z}$, i.e., $B\in \mathcal{L}(\mathcal{U},\mathcal{Z})$. The function $z$ is called the state of the system, with initial condition $z(0)=z_0$, while $u$ is the control input.

Let $G(s)$ be the transfer function of a linear and time-invariant infinite-dimensional system with input $u \in \mathcal{U}$ and output $y \in \mathcal{Y}$,  where $\mathcal U$ and $\mathcal Y$ are Hilbert spaces
\rev{ and
$$y(t) = C z(t) +D u(t) $$
where $C \in \mathcal{L} ( \mathcal{Z}, \mathcal{Y} ), $ $ D \in \mathcal{L} ( \mathcal{U}, \mathcal{Y} ). $  
Defining 
\begin{align*}
G(s) &= C (sI -A)^{-1} B + D ,\\
  \hat{y}(s) & = G(s)\hat{u}(s) \, . 
  \end{align*}
  }
    The space $H_{\infty}$ is defined as follows:
\begin{align*}
 H_{\infty} =
  \biggl\{ G: \mathbb{C}^+_0 \to \mathcal{L}(\mathcal{U},\mathcal{Y}) \, |& \, G \text{ is holomorphic and }\\
&\sup_{\text{Re} \, s>0} \| G(s)\| <\infty \biggr\}, 
\end{align*}
where $\mathbb{C}_0^+$ are all complex numbers with real part larger than zero. The norm on $H_{\infty}$, for $G\in H_{\infty}$, is
$$
\| G \|_{\infty} = \sup_{\text{ Re} \, s>0} \|G (s) \|,$$
and further it can be given by 
$$\| G \|_{\infty} = \sup_{\omega \in \mathbb{R}} \; \|G(j\omega)\|$$

\begin{defn}
 \label{defn-ext-stable}
 A system is {\em externally stable} or {\em $L_2$-stable} if for every input $u \in L_2 (0, \infty; \mathcal U)$, and zero initial condition,  the output $y \in L_2 (0,\infty ; \mathcal Y)$. 
  If a system is externally stable, the maximum ratio between the norm of the input and the norm of the output is called the {\em $L_2$-gain}. 
 \end{defn}
 
 The following results are well-known; see for instance \cite{curtain2012introduction}.
\begin{thm}
A system is {\em externally stable} if and only if its transfer function is in $H_\infty .$ 
\end{thm}

\begin{defn}
The system \eqref{sys:AB} is {\em stabilizable} if there exists $K \in \mathcal L (\mathcal Z , \mathcal U) $ such that $A+BK$ generates an exponentially stable semigroup.
\end{defn}

\begin{thm}
The system \eqref{sys:AB} with $y=z$ is externally stable if $A$ generates an exponentially stable semigroup. If the system is stabilizable, then external stability implies exponential stability of the semigroup generated by $A. $
\end{thm}

Finally, the following result on optimization in Hilbert space is included. 
\begin{thm}[{\cite[p.161, Th. 1]{luenberger1969}}] 
\label{minimumnorm}%
Let $\mathcal{X}$ and $\mathcal{Y}$ be Hilbert spaces and let $A\in \mathcal{L}(\mathcal{X},\mathcal{Y})$ with range closed in $\mathcal{Y}$. Then a vector $x \in \mathcal{X}$ of minimum norm satisfying $Ax=y$ for $y \in \mathcal{Y}$, if any exists, is given by $x=A^*z$ where $z$ is any solution of $AA^*z=y$. 
\end{thm}
The minimum norm solution to $Ax=y$ includes an operator $A^\dagger$ as defined in the following. 
\begin{defn}[{\cite[p.163]{luenberger1969}}]
Among all vectors $x_1\in\mathcal{X}$ satisfying $$\|Ax_1-y\|=\min_x \|Ax-y\|,$$ let $x_0$ be the unique vector of minimum norm. The pseudoinverse $A^\dagger$ of $A$ is the operator mapping $y$ into its corresponding $x_0$ as $y$ varies over $\mathcal{Y}$.
\end{defn}
\begin{rem}
Notice that $\tilde{y}=Ax_1$ is unique. The set of vectors $x_1$ that satisfy $\tilde{y}=Ax_1$ is a closed linear variety. Thus, it contains a unique $x_0$ of minimum norm. This makes $A^\dagger$ well defined. Furthermore, $A^\dagger$ is linear and bounded. If $A^*A$ is invertible, then a simple explicit formula can be given for $A^\dagger$ namely $A^\dagger=(A^*A)^{-1}A^*$. Similarly, if $AA^*$ is invertible, then $A^\dagger=A^*(AA^*)^{-1}$.  
\end{rem}

\section{Closed-Form $H_{\infty}$ Optimal State Feedback}  
\label{theorysection1}
 \label{secSF}
Consider the linear time-invariant system on the Hilbert space $\mathcal Z$ for $t \geq 0 ,$
\begin{align}
\frac{dz(t)}{dt} = Az(t) +Bu(t) +Hd(t), \quad z(0) = 0 , \label{sysSF}
\end{align}
where  $A : D(A) \subset \mathcal Z \to \mathcal Z$ is the generator of an infinitesimal semigroup on $\mathcal{Z}$ and  $B \in \mathcal L (\mathcal U, \mathcal Z)$ and $H \in \mathcal L (\mathcal D, \mathcal Z)$ where $\mathcal{U}$ and $\mathcal{D}$ are Hilbert spaces.
The system is affected by two exogenous signals, the control input $u(t)\in L_2 (0,\infty; \mathcal U)$ and  the disturbance $d(t)\in  L_2([0,\infty);\mathcal{D)} . $  The control signal is to be designed while the disturbance is unknown.

Consider, for a controller with transfer function $K(s),$ the state feedback  ${\hat{u}(s) = K(s)\hat{z}(s)}$. Furthermore, define the regulated output
\begin{equation*}
\hat \zeta = \begin{bmatrix} \hat z \\ R\hat u \end{bmatrix},
\end{equation*}
with $R \in \mathcal{L}(\mathcal{U}, \mathcal{R})$, where $\mathcal R$ is a Hilbert space, and such that $R^*R$ is coercive. The closed-loop transfer function from the disturbance $\hat{d}$ to $\hat{z}$, with the state feedback  ${\hat{u}(s) = K(s)\hat{z}(s)}$, is denoted $G_K$. That is, ${\hat \zeta (s) = G_K(s) \hat d(s)}.$ 
\rev{
\begin{defn}
The  {\em fixed attenuation $H_\infty$ control problem} with attenuation $\gamma$  is to construct a stabilizing controller with transfer function $K$ so that the closed loop system $G_{K}$   is externally stable and satisfies the bound 
\begin{equation}
 \| G_{K} \|_{\infty} ~ < ~ \gamma.
\label{eq:std0}
\end{equation}
If such a stabilizing controller  exists, the system is {\em stabilizable with attenuation $\gamma . $}
If $\gamma_{\text{opt}} = \inf \gamma $ over all $\gamma$ where the system is stabilizable with attenuation $\gamma ,$ then $\gamma_{\text{opt}}$ is the {\em optimal attenuation}. Moreover, 
$$\gamma_{\text{opt}} = \inf_K \; \|G_K\|_{\infty}$$
where the infimum is taken over the set of controllers $K$ for which the closed-loop system $G_K$ is stable. Finally, an optimal controller $K$, denoted $K_{\text{opt}}$, is one that fulfils $\|G_{K_{\text{opt}}}\|_{\infty} = \gamma_{\text{opt}}$ and for which $G_{K_{\text{opt}}}$ is stable. 
\end{defn}}

The following $H_{\infty}$ optimal problem is addressed; find a stabilizing controller $K$, if one  exists,  achieving optimal attenuation. 
The following theorem states a closed-form solution to this problem for a certain class of systems. 

\begin{thm}\label{mainTh} Consider \eqref{sysSF} and the definition of $G_K$ in the text after \eqref{sysSF}. 
 If $A$ is self-adjoint and strictly negative, then $$\inf_{K} \; \|G_K\|_{\infty}$$ is solved by the state feedback controller 
\begin{equation}
{K_{\text{opt}} = (R^*R)^{-1}B^*A^{-1}}
\label{eqn-Kopt}
\end{equation}
 and the minimum norm, or optimal attenuation,  is 
\begin{equation}
\gamma_{\text{opt}} = {{\|H^*(A^2+B(R^*R)^{-1}B^*)^{-1}H\|}}^{\frac{1}{2}}.
\label{eqn-gamma-opt}
\end{equation}
\end{thm}

\proof
The proof is divided into two parts. In the first part, it is shown that no controller can achieve a lower closed loop $L^2$-gain  than  \eqref{eqn-gamma-opt}. 
\rev{In the second part of the proof, it is shown that $G_K$ with $K = K_{\text{opt}}$ is stable and that $\|G_{K_{\text{opt}}}\|_{\infty}$ is equal to  \eqref{eqn-gamma-opt}.}

(1)
By the definition of the $H_{\infty}$ norm $$\|G_K\|_{\infty} = \sup_{\omega \in \mathbb{R}} \|G_K(j\omega)\|.$$ 
It follows that 
$$\|G_K\|_{\infty} \geq \|G_K(0)\|,$$
A  $K\in \mathcal{L}(\mathcal{Z},\mathcal{U})$ that minimizes $\|G_K(0)\|$ will now be calculated. 
Note that 
\begin{equation*} 
\|G_K(0)\| = \sup_{d \in \mathcal{D}, \; \|d\| = 1} \; \|G_K(0)d\|. 
\end{equation*}
This can be  written equivalently  as 
\begin{align}
\sup_{d \in \mathcal{D}, \,\rev{\| d\|= 1}}& \quad \|\zeta\|\\ 
\text{subject to}& \quad 0 =Az + Bu+Hd, \\ 
& \quad \zeta = \begin{bmatrix} z \\ R u \end{bmatrix}, \\ 
& \quad u = Kz. \label{controller}
\end{align}
Now, without the constraint \eqref{controller} the problem becomes
\begin{equation}\label{relaxation}
\begin{aligned}
\sup_{d \in \mathcal{D}, \,\rev{\| d\|= 1}}& \quad \| \zeta\|\\ 
\text{subject to}& \quad 0 =\begin{bmatrix}A \quad& B(R^*R)^{-1}R^*\end{bmatrix}{\zeta} +H{d}.
\end{aligned}
\end{equation}
Further, the property of strict negativity implies that $A$ has a bounded inverse. 
Define $$T \coloneqq -\begin{bmatrix}I \quad& A^{-1}B(R^*R)^{-1}R^*\end{bmatrix}.$$
For $z$ in the range of $A$ 
$$\|A^{-1}z\|\leq \frac{1}{m}\|z\|$$
is implied by the strict negativity of $A$, that is $\langle -Az,z\rangle \geq m \|z\|^2$ for some $m>0$.  It follows that $A^{-1}$ is bounded on the range of $A$. Now, let $x$ be such that $\langle Az,x\rangle = 0$ for all $z\in D(A)$. The definition of an adjoint operator implies that $A^*x = 0$ for $x \in D(A^*)$. Since $A$ is self-adjoint, we conclude that $Ax=A^*x= 0$. The negativity of $A$ shows that this can only happen if $x =0$, so the range of $A$ is dense in $\mathcal{Z}$. Thus, $A^{-1} \in \mathcal{L}(\mathcal{Z})$. Therefore, $T$ is a bounded linear operator on $\mathcal Z .$

Now, consider 
$$TT^* = I+A^{-1}B(R^*R)^{-1}B^*A^{-1}.$$
Notice that the operator $TT^*$ is strictly positive and thus invertible, with bounded inverse. 
Thus, by Theorem~\ref{minimumnorm} 
 the constraint in  \eqref{relaxation} can be rewritten as ${\zeta = T^{\dagger}A^{-1}H d}$
where 
\begin{align*}
T^{\dagger} &\coloneqq T^*(TT^*)^{-1} =\\& -\begin{bmatrix} I\\R(R^*R)^{-1}B^*A^{-1}\end{bmatrix} \left(I+A^{-1}B(R^*R)^{-1}B^*A^{-1}\right)^{-1}.
\end{align*}
Moreover, 
\begin{align*} 
\|G_K(0)\|  &=\sup_{d\in \mathcal{D},\, \rev{\| d\|= 1}} \;\|T^{\dagger}A^{-1}H d\|  \\ & = \|T^{\dagger}A^{-1}H\|  \\ & = \|H^*A^{-1}(T^{\dagger})^*T^{\dagger}A^{-1}H\|^{\frac{1}{2}} \\
& = \|H^*\left(A^2+B(R^*R)^{-1}B^*\right)^{-1}H\|^{\frac{1}{2}} .
\end{align*}
This completes the first part of the proof. 

(2) Define $K_{\text{opt}}$  as in \eqref{eqn-Kopt}  and let $T(t)$ be the strongly continuous semigroup on $\mathcal{Z}$ generated by $A+BK_{\text{opt}}.$ 
Setting $P=-A^{-1} ,$
\begin{align*}
&\langle (A+B(R^*R)^{-1}B^{*}A^{-1})z,Pz \rangle \\ &\qquad\qquad \qquad+\langle Pz,(A+B(R^*R)^{-1}B^*A^{-1})z \rangle \\&= 
\langle (-2I-2A^{-1}B(R^*R)^{-1}B^{*}A^{-1})z,z \rangle\\
&\leq -\langle z,z\rangle . 
\end{align*}
Theorem \ref{lemLyap} then implies that $T(t)$ is exponentially stable. In other words, $K_{\text{opt}}$ is stabilizing. Finally, since 
\begin{multline*}
G_{K_\text{opt}}(s) =\\ \begin{bmatrix}I\\ R(R^*R)^{-1}B^*A^{-1} \end{bmatrix}(sI-A-B(R^*R)^{-1}B^*A^{-1})^{-1}H,
\end{multline*}
it will be shown that $$ \|G_{K_{\text{opt}}}\|_{\infty}  \leq \|H^*\left(A^2+B(R^*R)^{-1}B^*\right)^{-1}H\|^{\frac{1}{2}}.$$
Notice that the inequality above is equivalent to 
\begin{multline*}
\langle d,G_{K_{\text{opt}}}( \jmath \omega)^*G_{K_{\text{opt}}}( \jmath \omega)d\rangle \\ \leq \langle  d, H^*\left(A^2+B(R^*R)^{-1}B^*\right)^{-1}Hd \rangle
\end{multline*}
for all $d\in \mathcal D$. Define $M \coloneqq A^2+B(R^*R)^{-1}B^*$. The inequality  can now be stated as 
\begin{multline*}
\langle d,H^*\left( (\jmath \omega A-M)M^{-1} (- \jmath \omega A-M) \right)^{-1}Hd\rangle \\\leq \langle  d, H^*M^{-1}Hd \rangle,
\end{multline*}
and further simplified to
\begin{equation*}
\langle d,H^*\left(\omega^2 AM^{-1}A+M \right )^{-1}Hd\rangle \\ 
\leq \langle  d, H^*M^{-1}Hd \rangle.
\end{equation*}
From the latter representation of the inequality it is evident that it holds for all $\omega \in \mathbb{R}$. This completes the proof.

The result in Theorem~\ref{mainTh} is a rarity in $H_{\infty}$ synthesis for infinite-dimensional systems as closed-form expressions are generally hard to obtain. Furthermore, as for finite-dimensional systems,  it is in general necessary to compute  a series  of suboptimal controllers to find the controller achieving optimal attenuation.

\section{Closed-Form $H_{\infty}$ Optimal Filtering}  
\label{theorysection2}
Consider now the observation system
\begin{equation} \label{sysSE}
\begin{aligned}
\frac{dz(t)}{dt} &= Az(t)+w(t), \quad z (0) = 0 \\
q(t) &= Qz(t),\\
y(t) &= Cz(t)+\rev{S}v(t) . 
\end{aligned}
\end{equation}
Again, $A$ generates a strongly continuous semigroup on $\mathcal Z , $ $Q \in \mathcal L ( \mathcal Z , \mathcal{Q} ) $ and $C \in  \mathcal L ( \mathcal Z , \mathcal{Y} )$. In this class, $y(t)\in \mathcal{Y}$ are the available measurements while $w(t)\in \mathcal{Z}$ and \rev{$v(t)\in \mathcal{V}$, where $\mathcal{V}$ is a Hilbert space}, are process disturbance and measurement noise, respectively. \rev{Finally, $S \in \mathcal L ( \mathcal V , \mathcal{Y} ) $.}

The objective is to design an observer that optimally estimates  $q(t)$, which equals $z(t)$ if $Q=I,$  given disturbance and noise in $L_2([0,\infty);\mathcal{Z})$ and $L_2([0,\infty);\mathcal V)$. That is, the estimator  state  $\nu $ minimizes the $H_{\infty}$ norm of the transfer function from the disturbance and noise to the estimation error $Q z - Q \nu . $ Such an observer is given in closed form in the following theorem.   

\begin{thm} \label{mainTh2}
Consider operators $A: D(A) \rightarrow \mathcal{Z}$, ${Q \in \mathcal{L}\, (\mathcal{Z},\mathcal{Q})}$, ${C \in \mathcal{L}\, (\mathcal{Z},\mathcal{Y})}$ \rev{and $S \in \mathcal L ( \mathcal V , \mathcal{Y} ) $}, where $\mathcal Z$, $\mathcal Q$, $\mathcal{Y}$ and \rev{$\mathcal{V}$} are Hilbert spaces. Consider the class of observers
\begin{align*}
\frac{d\nu(t)}{dt} &= A\nu(t)+L(C\nu(t)-y(t)), \quad \nu (0) = 0, 
\end{align*}
where $L\in \mathcal{L}(\mathcal Y,\mathcal Z)$. If $A$ is self-adjoint and strictly negative, then the estimation error
$$\sup_{\|[w^* \; v^* ]\|_2\rev{= }1} \|Qz-Q\nu\|_2,$$ 
is minimized by 
$$L_{\text{opt}} = A^{-1}C^*(SS^*)^{-1}.$$
The minimal estimation error  is $$\|Q(A^2+C^*\rev{(SS^*)^{-1}}C)^{-1}Q^*\|^{\frac{1}{2}}.$$
\end{thm}
\proof
If the disturbance and noise signals are stacked as 
$$ d (t) =\begin{bmatrix}w (t) \\ v (t)\end{bmatrix},$$
then the dynamics of the estimation error ${e(t) = z(t)-\nu(t)}$ and alternative output estimation error ${\zeta(t) = q(t)-Q\nu(t)}$ can be written as 
\begin{equation}\label{closedT}
\begin{aligned}
\frac{de(t)}{dt} &= (A+LC)e(t) +\begin{bmatrix} I \quad& L\rev{S} \end{bmatrix} d(t) ,\\
\zeta(t) &= Qe(t).
\end{aligned}
\end{equation}
Denote the transfer function from $\hat d$ to $\hat \zeta$ in \eqref{closedT}, with a specific $L$, by $G_L$. That is, $\hat{\zeta}(s)=G_L(s)\hat{d}(s)$. Note that
$$\|G_L\|_{\infty} = \sup_{\|\hat{d}\|_2\rev{=} 1} \|\hat{\zeta}\|_2 = \sup_{\|[\hat{w}^* \; \hat{v}^* ]\|_2\rev{=} 1} \|Q\hat{z}-Q\hat{\nu}\|_2 .$$
This is exactly the problem considered in the theorem statement. Consider the dual system of $G_L$, i.e., its adjoint $G_L^*(s)$. It holds that $\|G_L^*\|_{\infty} = \|G_L\|_{\infty}$. Given ${\hat{\zeta}(s) = G_L^*(s)\hat{d}(s)}$ with state space representation
 \begin{equation}\label{dual}
\begin{aligned}
\frac{dz(t)}{dt} &= (A+C^*L^*)z(t) +Q^*d(t) ,\\
\zeta(t) &= \begin{bmatrix}I\\\rev{S^*}L^* \end{bmatrix}z(t).
\end{aligned}
\end{equation}
This system is exactly of the form of the closed-loop system treated in Theorem~\ref{mainTh}, however, with $K = L^*$, $B=C^*$, $H = Q^*$ \rev{and $R= S^*$}. Thus, it follows from Theorem~\ref{mainTh}, and the fact that $\|G_L\|_{\infty}=\|G_L^*\|_{\infty}$, that $L^* = (SS^*)^{-1}CA^{-1}$ is optimal, i.e., $L = A^{-1}C^*(SS^*)^{-1}= L_{\text{opt}}$. \rev{Finally, 
as $\|G_L\|_{\infty} = \|G_L^*\|_{\infty}$ and the minimal value of $\|G_L^*\|_{\infty}$ is 
$$\|Q(A^2+C^*(SS^*)^{-1}C)^{-1}Q^*\|^{\frac{1}{2}},$$
by Theorem~\ref{mainTh}, it holds that the minimal value of $\|G_L\|_{\infty}$  is also 
$$\|Q(A^2+C^*\rev{(SS^*)^{-1}}C)^{-1}Q^*\|^{\frac{1}{2}}.$$}

\section{Temperature Regulation with Variable Conductivity} 
\label{variableconductivity}

\begin{figure}
\begin{center}
\includegraphics{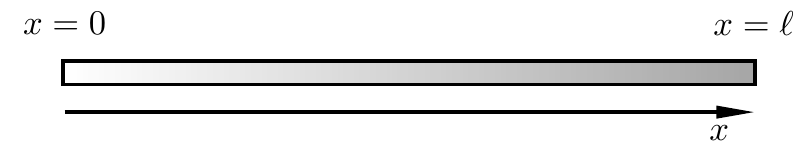}
\caption{Rod of length $\ell$ with one-dimensional spatial coordinate $x$ and conductivity $k(x)=x^2+1$.}
\label{fig:heatrod}
\end{center}
\end{figure}

Consider the one-dimensional rod of length $\ell$ in \figref{fig:heatrod}. The material of the rod has variable conductivity as described by $k(x)=x^2+1$, where $x$ is the spatial coordinate. The heat propagation in the rod is governed by 
\begin{equation} \label{heatEqSF1}%
\begin{aligned}
\frac{\partial z}{\partial t}(x,t) = \frac{\partial}{\partial x}\left(k(x)\frac{\partial}{\partial x}\right)z(x,t) +Bu(t)+d(t),\\\quad 0 < x < \ell, \; t \geq 0,
\end{aligned}
\end{equation}
where $z(x,t)$ is the temperature at time $t$ in position $x$. The temperature profile can be regulated by the control input $u(t)$ and is disturbed by the signal ${d(t) \in L_2([0,\infty);\mathbb{R})}$.  The temperature is kept at zero at the end points. That is, Dirichlet boundary conditions $z(0,t) = 0$, $z(\ell,t) = 0,$ hold. Furthermore, $Bu(t) = u(t)$ for  $0<x<\ell$. 

Define the operator $A$ as
\begin{align} \label{Aop}
A = \frac{\partial}{\partial x}\left(k(x)\frac{\partial}{\partial x}\right),
\end{align}
with domain
\begin{equation*}
\begin{aligned}
D(A) =  \biggl\{ z \in L_2(0,\ell) \, \biggl |  &\,z \text{ and }\frac{dz}{dx} \text{ locally absolutely} \\ &\text{continuous in } [0,\ell ], \\ &\frac{d^2 z}{dx^2} \in L_2(0,\ell) \\ &\text{ with } z(0) = 0, \; z(\ell) = 0 \biggr \}.
\end{aligned}
\end{equation*}
Then, we can write \eqref{heatEqSF1} on abstract form
\begin{align*}
\dot{z}(t) = Az(t)+Bu(t)+d(t). 
\end{align*}
The operator $A$ is self-adjoint and strictly negative on $D(A)$, see e.g. \cite{curtain2012introduction}. Hence, by Lemma~\ref{lem1}, $A$ is the infinitesimal generator of a strongly continuous semigroup. The system is thus of the form considered in Theorem~\ref{mainTh}, \rev{with $R = I$,} and an $H_{\infty}$ optimal state feedback is 
 \begin{align*}
 u(t) &= B^*A^{-1}z(t)  \\& = \int_{0}^{\ell}\int_0^\ell G(x,s)z(s,t) ds \, dx 
 \end{align*}
where
\begin{align*}
G(x,s) = \begin{cases} \left(\frac{\arctan(s)}{\arctan (\ell)}-1 \right)\arctan (x)& \text{if } \, 0 < x < s,  \\[10pt] \left(\frac{\arctan(x)}{\arctan (\ell)}-1 \right)\arctan (s) &\text{if } s<x <  \ell. \end{cases} \, 
\end{align*}
The function $G(x,s)$ is the Green's function of the linear differential operator $A$ in \eqref{Aop}, when ${k(x)=x^2+1}$. See \cite{zaitsev2002handbook} for more details on Green's functions for Sturm-Liouville problems. 
Defining
$$ g(s) = \frac{1}{2}\left(\ln(s^2+1)-\ln(\ell^2+1)\frac{\arctan(s)}{\arctan(\ell)} \right) ,  $$
the optimal control is
\begin{equation}\label{gs}
    u(t) = \int_0^\ell  g(s)  z(s,t) \,ds.
\end{equation}
 If the conductivity was constant; that is, $k(x) = 1$, $g(s)$ would simply be ${g(s) = s(s-\ell)/2}$. See \cite{lidstrom2016h} for an example with constant conductivity.
 
\begin{figure}
\begin{center}
\includegraphics{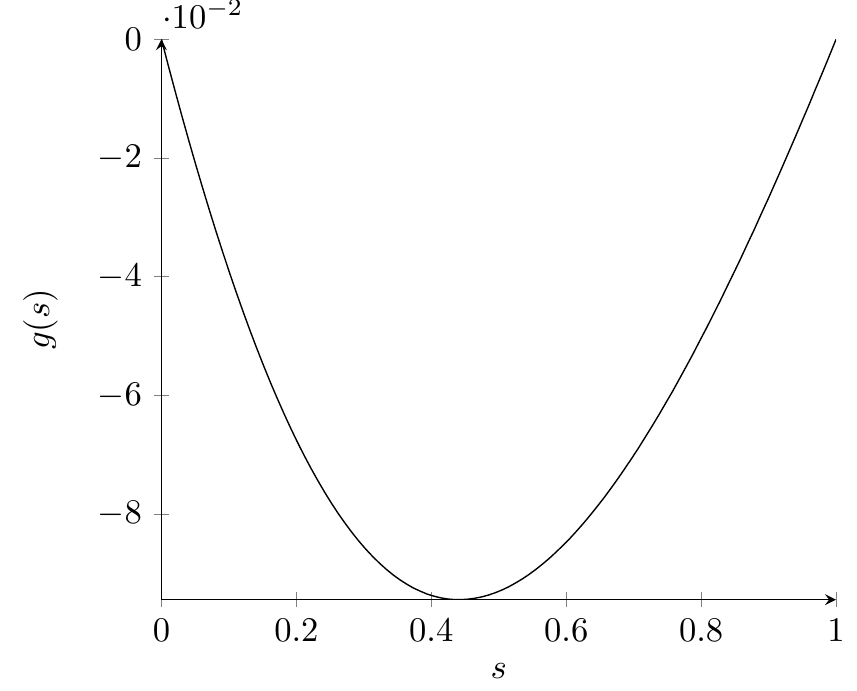}
\caption{Plot of feedback weight gain $g(s)$ in \eqref{gs} when $\ell  = 1$.}
\label{fig:plotg}
\end{center}
\end{figure}
The feedback gain $g(s)$ in \eqref{gs} is plotted in \figref{fig:plotg} for the case with $\ell = 1$. It determines the scalar signal for controlling the temperature profile, as a compromise between the deviation in temperature from zero and the cost for changing the temperature by means of the control signal, weighted by the conductivity. Note that it is not symmetric across $s$ due to the conductivity $k(x)=x^2+1$ not being symmetric along the rod. More specifically, the conductivity increases towards the right-hand side of the rod in Figure~\ref{fig:heatrod}. If we compare the absolute value of the gain $g(s)$ at $s= 0.5-a$ and $s = 0.5+a$, for any $0<a<0.5$, it is lower at the latter point. This is due to that the material further to the right is better at transporting heat.


\section{Estimation of Plate Temperature} 
\begin{figure}
\begin{center}
\includegraphics{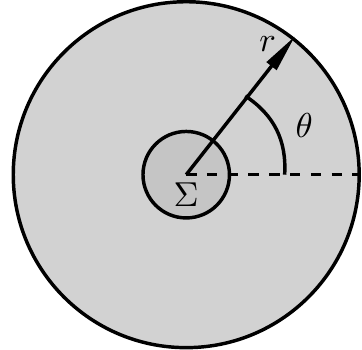}
\caption{Circular plate with radius $1$ and spatial coordinate $(r,\theta)$.}
\label{fig:heatplate}
\end{center}
\end{figure}

Consider the following partial differential equation that models heat propagation in a circular plate, with radius of length $1$, given a disturbance ${w(t) \in L_2([0,\infty);\mathbb{R})}$,
\begin{equation}\label{heatEqEst}
\begin{aligned}
\frac{\partial z}{\partial t}(r,\theta,t) = \left(\frac{\partial ^2 }{\partial r^2}+\frac{1}{r}\frac{\partial }{\partial r}+\frac{1}{r^2}\frac{\partial ^2 }{\partial \theta^2}\right)z(r,\theta,t) +w(t), \\ 0 \leq r < 1,\quad  0 \leq \theta \leq 2 \pi, \quad  t \geq 0. 
\end{aligned}
\end{equation}
Here, $z(r,\theta,t)$ is the temperature at time $t$ and position $(r,\theta)$, see \figref{fig:heatplate} for a depiction of the plate. The plate area, not including the boundary, is denoted $\Omega$. Moreover, the temperature is kept at zero on the boundary, i.e.,
  $$z(1,\theta,t) = 0 , \quad 0\leq \theta \leq 2\pi .$$
  
 The average temperature can be measured over a circular region of radius $1/4$, with its centre aligned with that of the plate. The region is denoted $\Sigma$, see \figref{fig:heatplate}.  Define
  \begin{equation} \label{Cop}
Cz =  \int_{\Sigma} z(r,\theta,t) r\,dr\,d\theta.
\end{equation}
The measurement can then be constructed as 
\begin{equation}
y(t) = Cz(r,\theta,t)+ v(t), \label{measurements}
\end{equation}
where $v(t) \in L_2([0,\infty);\mathbb{R})$ accounts for sensor errors.  
Theorem~\ref{mainTh2} will now be used to obtain an observer  that optimally estimates the temperature profile $z$ based on the measurements $y$. Note that $Q = I$ and $S=I$ in this example. 

First, define 
\begin{align} \label{Aop2}
A = \frac{\partial ^2 }{\partial r^2}+\frac{1}{r}\frac{\partial }{\partial r}+\frac{1}{r^2}\frac{\partial ^2 }{\partial \theta^2}
\end{align}
with domain
\begin{equation*}
\begin{aligned}
D(A) &=\\ 
 \biggl \{& z \in L_2\left(\Omega \right) \, \biggr \rvert z\text{ locally absolutely continuous in } \Omega, \\ & \text{given } x = r\cos(\theta), \; y = r\sin(\theta), \;  \frac{\partial^k z}{\partial x^k}, \, \frac{\partial^k z}{\partial y^k} \text{ for } \\ &k = 1,2,\text{ locally absolutely continuous in }\Omega. \; \\&Az \in  L_2\left (\Omega\right) \text{ with }
z(1,\theta,t) = 0 \text{  for  } 0\leq \theta \leq 2 \pi \biggr \}  .
\end{aligned}
\end{equation*}
The operator $A$  is  self-adjoint and strictly negative. Hence, it generates an exponentially stable strongly continuous semigroup on $L_2\left (\Omega \right).$ Moreover, the state $z$ evolves on the space $L_2\left(\Omega \right)$ and we can write \eqref{heatEqEst} and \eqref{measurements} abstractly as 
\begin{equation}\label{Exest}
\begin{aligned} 
\dot{z}(t) &= Az(t)+w(t),\\
y(t) &= Cz(t)+v(t).
\end{aligned}
\end{equation}

The system \eqref{Exest} fulfils the requirements of Theorem~\ref{mainTh2}.
Thus, the optimal observer is 
\begin{align*}
\frac{d\nu(t)}{dt} &= A\nu(t)+L_{\text{opt}} (C\nu(t)-y(t)), \quad \nu(0)=0,
\end{align*}
where $L_{\text{opt}} = A^{-1}C^*$. As in the previous subsection, we can obtain an expression for $$u(r,\theta,t)\coloneqq L_{\text opt}y(t).$$
Defining
 \begin{align*}
G(r,\theta,r_0,\theta_0)= \frac{1}{4\pi}\ln \frac{r^2+r_0^2-2rr_0\cos(\theta-\theta_0)}{r^2r_0^2+1-2rr_0\cos(\theta-\theta_0)},
 \end{align*}
 then 
  $$
 L_{\text{opt}}  y(t) = \int_{0}^{2\pi}\int_{0}^{\frac{1}{4}} G(r,\theta,r_0,\theta_0) r_0dr_0d\theta_0  \cdot y(t).$$
 Note that the problem is rotationally invariant. Furthermore, $u(r,\theta,t)$ can also be determined as the solution to the Laplace equation on $\Omega$, subject to the given boundary condition, where there is a uniformly distributed source in $\Sigma$ (i.e., it is a Poisson's equation on $\Sigma$). 
 
 
  
 \section{Implications for General Purpose Algorithms} \label{numericalsection}
For finite-dimensional systems, the $H_{\infty}$ optimal synthesis is generally performed though iteratively solving a series of AREs, as mentioned in the introduction.  The method in \cite{arnold1984generalized} is one example of an algorithm for solving AREs. However, it works poorly when the order of the system is large. Other methods for solving AREs are  the matrix sign function method \cite{byers1987solving} and the game-theoretic method developed in \cite{lanzon2008computing}. Even though synthesis techniques with AREs have existed for decades, there is no generally accepted algorithm for systems of large order. The sign-indefiniteness of the quadratic term in the $H_{\infty}$-ARE and the need for an iterative procedure to find the optimal attenuation complicate computation. 

The results stated in Theorem~\ref{mainTh} and \ref{mainTh2} can be used in benchmarking of general purpose algorithms for $H_{\infty}$ control. Now, the computational time of computing the closed-form controller in Theorem~\ref{mainTh} will be compared with the computational time of a general purpose algorithm for $H_{\infty}$ synthesis. More specifically, it will be compared with the method for solving AREs described in \cite{arnold1984generalized}. 

The comparison is performed using Matlab~R2012a \cite{MATLAB} on an Intel Core i7-3770K CPU with 3.5 GHz $\times 8$ and 31.3 GiB memory. The methods in \cite{arnold1984generalized} is implemented in Matlab as the function CARE. An optimal state feedback controller will be obtained through bisection, where the $H_2$-control solution is used as the upper bound and zero is used as the lower bound. 

\begin{figure} 
\centering
\includegraphics[scale=1]{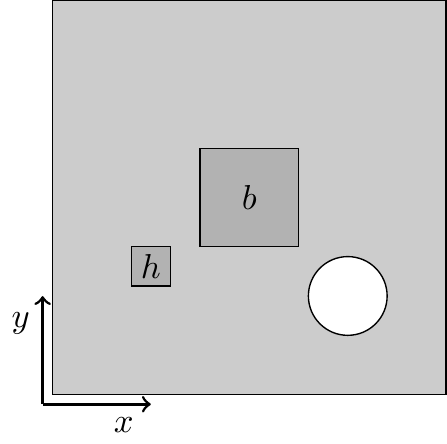}
\caption{Irregular geometry considered in algorithm comparison. \label{planes}}
\end{figure}

The comparison will be performed on an example system of heat diffusion on a two-dimensional irregular geometry. The geometry is a plane of $4 \times 4$ units, with a circle of radius $0.4$ units at $(3,1)$ removed. The lower left corner of the plane is the origin. See Figure~\ref{planes} for a depiction of the plane. Inspiration for this example is taken from \cite{kasinathan2014solution}. 

The heat distribution at position $(x,y)$ at time $t$ is denoted $z(x,y,t)$ and governed by
\begin{equation}
\begin{aligned}
        \frac{\partial z}{ \partial t}(x,y,t) &= \frac{\partial^2 z}{\partial x^2}+\frac{\partial^2 z}{\partial y^2} +b(x,y)u(t)+h(x,y)d(t),
\end{aligned} \label{heatEq}
\end{equation}
with zero Dirichlet boundary condition. Furthermore, the functions $b(x,y)$ and $h(x,y)$ are defined as 
\begin{align*}
    b(x,y)= \begin{cases} 1 & (x,y) \in \{\textrm{ square at } (2,2) \textrm{, side } 1 \},\\ 0 & \text{otherwise.}\end{cases}
\end{align*}
and
\begin{align*}
    h{(x,y)}= \begin{cases}1 & (x,y) \in \{\textrm{ square at } (1,1.3)\textrm{, side } 0.4 \},\\ 0 & \text{otherwise.}\end{cases}
\end{align*}
\rev{The squares are centered at the given positions and have the specified side-length, e.g., a square at (2,2), side 1 is a square centered at (2,2) with side-length 1, see Figure~\ref{planes}. }

The infinite-dimensional system \eqref{heatEq} is approximated by a finite-element method with linear splines. The resulting finite-dimensional system, for a given approximation order, is on the following form 
\begin{equation} \label{genSys}
    \begin{aligned}
    \dot{z} &= Az+Bu+Hd, \\ 
    \zeta &= \begin{bmatrix}z \\ u \end{bmatrix},
\end{aligned}
\end{equation}
where $z(t)\in \mathbb{R}^n$ is now the finite-dimensional state of the system and $n$ is the order of the system. Furthermore, $A$, $B$, $H$ are matrices of appropriate dimensions \rev{with $A = A^T$, $A$ negative definite. Thus, $A$ is self-adjoint and strictly negative and  Theorem~\ref{mainTh} applies. }

 \begin{figure}
\begin{center}
\includegraphics[scale=1]{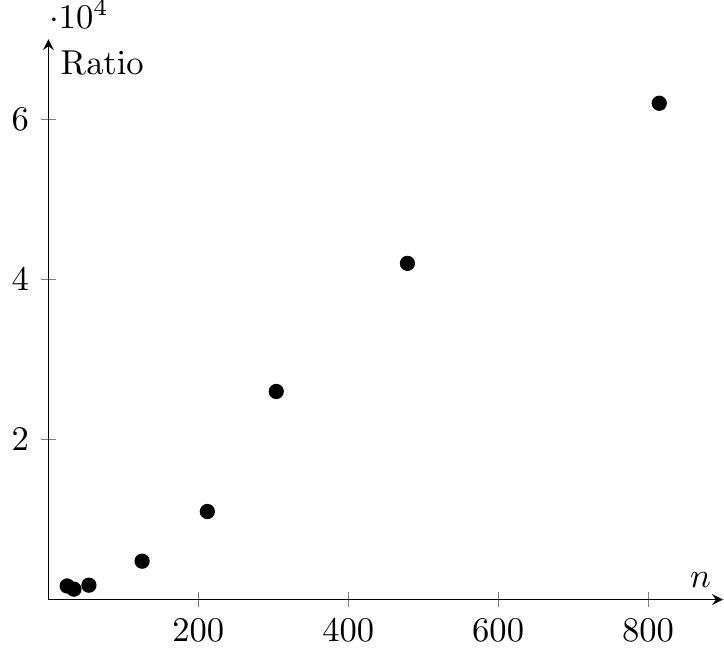}
\caption{Ratio between the computational time of  the method  for calculating an optimal $H_\infty$-controller proposed in this paper  and the standard iterative approach, for various approximation orders $n$. For $n = 815$, the computational time is 19.4 ms for the proposed method and 20.2 minutes for the standard method.}
\label{fig:plotAlgComp}
\end{center}
\end{figure}

The \rev{ two methods for computation of the optimal attenuation and control are:}
\begin{enumerate}[I]
    \item Compute the optimal controller $K_{\textrm{opt}} =B^TA^{-1}$, 
    \item Compute an optimal state feedback $K$ through solution of a series of AREs and a  bisection algorithm. \rev{The function CARE in Matlab, as described in \cite{arnold1984generalized}, was used to solve the AREs.}
\end{enumerate}

The result of the comparison is given in \figref{fig:plotAlgComp}. The new method I is at least $10^3$ times faster than the standard iterative approach of method II. This improvement in computational speed has implications for lowering the computational time for large-order systems. \rev{The code for the comparison can be found at \url{https://gitlab.control.lth.se/carolina/code-hs}}.

\section{Conclusions and Directions for Further Research}
Simple expressions for $H_{\infty}$ state feedback and estimation  in closed form were obtained. They are optimal for infinite-dimensional systems with self-adjoint and strictly negative generator. Diffusion equations are an important example in this class. The results were illustrated with several analytical examples and also compared with a standard numerical algorithm. The results can also be used for testing and benchmarking of general purpose algorithms for $H_{\infty}$ synthesis. \rev{The proposed state feedback and estimation operators may be approximated in finite-dimensions before implementation. A procedure for how to do this as well as guarantees for the approximated gains to be admissible are not included in this paper. Instead, we refer to \cite{Morris-DPSbook} for  details. }

Future work includes investigation of how the closed-form expressions can be used to improve the performance of $H_{\infty}$ synthesis algorithms for general systems. 
Theorems ~\ref{mainTh} and \ref{mainTh2} may be useful for more general systems. \rev{The frequency response of the systems considered in this paper have maximum gain at zero frequency. This is used in the proofs of Theorem~\ref{mainTh} and \ref{mainTh2}, respectively. The optimal state feedback is constructed using the controller that  minimizes the gain at zero frequency. A similar  statement holds for the observer.
This  reveals a possible  generalization. In the case of state feedback the generalization is $u = B^*A^{-*}z$, where the assumption on self-adjointness of $A$ is disregarded. The criteria for this more general feedback gain to be optimal are that it is stabilizing and that the gain of the closed-loop system is obtained at zero frequency. 
The class of systems to which this control law is applicable is not restricted to those with self-adjoint $A$ \cite{bergeling2019closed}.}
This lower bound can be used to narrow the bisection range for more general systems.  The general purpose method for $H_{\infty}$ synthesis described in \cite{lanzon2008computing} constructs a unique monotonic matrix sequence, initialized with the zero solution, which converges to the stabilizing solution of the ARE.  The control law $u = B^*A^{-*}z$ may be a good initial iterate for systems \rev{where $\|A-A^*\|$ is small, i.e. $A$ is close to self-adjoint. However, this needs to be investigated further.}

\rev{Moreover, in control of diffusion and other  systems with partial differential equation models, the location of the actuator is generally a design variable. This is also true for the placement of the sensor. The location should be chosen so as to optimize the considered performance objective. Preliminary investigations suggest that the closed-form expression for the optimal attenuation could be utilized for these design problems. In \cite{morris2011linear}, the design of linear-quadratic optimal actuator placement is made on approximations of the PDEs and conditions for convergence are stated. Furthermore, an algorithm for computing the optimal location is given in \cite{darivandi2013algorithm}. Optimal actuator location measured by the $H_{\infty}$-norm is treated in \cite{kasinathan2013h}. However, the methods are computationally intensive. In particular, they involve many computations of the approximating $H_{\infty}$ synthesis problem. The closed form expressions in Theorem~\ref{mainTh} and \ref{mainTh2} could be used for efficient computation of the optimal actuator or sensor placement, respectively.  Moreover, the expression for the optimal $H_\infty$-norm can be minimized to find the best location. Similarly, one can parametrize the minimal error for the sensor location. Although the results are restricted to a certain class of systems, as for controller synthesis, this approach can be used in benchmarking of general purpose algorithms for optimal actuator and sensor location.}

\begin{ack}                               
Carolina Bergeling and Anders Rantzer are members of the LCCC Linnaeus Center and the ELLIIT Excellence Center at Lund University. This research was supported by the Swedish Research Council through LCCC and by the Swedish Foundation for Strategic Research. 
Furthermore, this paper is a result of  research initiated during the thematic year on Control Theory and its Applications 2015-2016, at the Institute for Mathematics and its Applications, University of Minnesota, Minneapolis, USA, in which all three authors took part.
\end{ack}

\bibliographystyle{plain}        
\bibliography{ref}           



\end{document}